\documentclass[12pt,a4paper]{article}
\usepackage{graphicx}
\usepackage[cp1251]{inputenc}
\usepackage[english]{babel}
\usepackage{amssymb, amsfonts, amsmath}

\begin{document}

\title{Classification of Einstein spaces with Stackel metric of type (3.0)}

\author{V. V.  Obukhov}

\maketitle

%
%
%


\quad

Tomsk  State  Pedagogical University, Institute of Scietific Research and
Development,, 60 Kievskaya St., Tomsk 634041, Russia; obukhov@tspu.edu.ru\\

Tomsk State University of Control Systems and Radio Electronics,Laboratory for Theoretical Cosmology, International Center of Gravity and Cosmos, 36, Lenin Avenue, Tomsk 634050, Russia

\quad


\section{Introduction}

In the middle of the last century, the problem of exact integration of Einstein equations attracted increased attention of gravitationalists, due to which a large number of exact solutions of these equations were found. Many of them played an important role in understanding the nature of gravitation and in the development of mathematical methods of the General Theory of Relativity. At present, the interest in finding new exact solutions has waned. On the one hand, the physical meaning of most known exact solutions is still unclear and requires additional research. On the other hand, the set of known exact solutions is not fully ordered, since a complete classification of known solutions has not been carried out, in order to separate from their set of all non-equivalent with respect to the transformation of coordinates solutions. Therefore, the classification of exact solutions of the Einstein equations, which allows us to identify and organize these solutions according to physically and mathematically understandable conditions, is still of considerable interest. From the methods of classification of exact solutions of Einstein equations noted in the book [1], we focused on the classification of Einstein spaces on the basis of their spatial symmetry. The set of such spaces includes Petrov spaces and Stackel spaces.
The problem of classification of Riemannian by groups of motions (Petrov spaces classification) has been solved in the works of A.Z. Petrov (see, for example, \cite{2}). In canonical coordinate systems all corresponding metrics has been found, and a partial classification of Ricci-flat spaces and Einstein spaces has been carried out.
In the Stackel spaces the privileged coordinate systems exist, in which the geodesic equations can be integrated by the method of complete separation of variables in the free Hamilton-Jacobi equation. On some subsets of the set of Stackel spaces the Lorentz, Klein-Gordon-Fock, Dirac-Fock and Weyl equations also can be integration by the method of complete separation of variables. Classifications of Stackel spaces are based on theorems on the separation of variables in the equations of motion of a test particle proved by V.N. Shapovalov (\cite{3}, \cite{4}, \cite{5}, see also \cite{5a}, \cite{5b}, \cite{5c}, \cite{5e}) . In these theorems, in particular, it is established that a necessary condition for complete separation of variables is the existence of a complete set in space consisting of mutually commuting Killing vector and tensor fields. It is proved that such properties are possessed by  Stackel spaces admitting privileged coordinate systems in which complete separation of variables is possible. In the case of a charged particle for the Hamilton-Jacobi equation, as well as for quantum mechanical equations, additional restrictions ( \cite{9}, \cite{10}, \cite{11}, \cite{12} ) are imposed on the complete set.
Privileged coordinate systems can be related by admissible coordinate transformations. Therefore, the corresponding Stackel spaces form equivalence classes. The task of classification is to enumerate all these classes. Note that the first works on the classification of the Stackel metrics (\cite{9} -\cite{12}) were devoted to flat spaces. As a result of such classification, all classes of privileged coordinate systems and complete sets of integrals of motion for the Hamilton-Jacobi, Klein-Gordon and Dirac equations have been listed.
The importance of realizing the classification of metrics of spaces by groups and symmetry algebras is evidenced by the fact that all the most interesting from the physical point of view of exact solutions of Einstein's equations, such as the static spherically symmetric Schwarzschild (\cite{13},\cite{14}), Kerr (\cite{15}), Newman-Unti-Tamburino (\cite{17}), Friedmann (\cite{18}), and some others (for example de Sitter cosmological solution), belong
to the number of Stackel spaces (subset of the set of spaces with symmetries). The interest in these spaces has not waned until now. As an example, one can point to the works in the field of cosmology, including the modified theory of gravitation {\cite{20}, \cite{21}, \cite{22}, \cite{22a}, \cite{22b}, \cite{23}, \cite{24}; and in other directions of research \cite{25}, \cite{25a}, \cite{25b}, \cite{25c}.

\noindent
Note that it is very difficult to solve the classification problem by investigating already known exact solutions of the Einstein equations due to the following circumstances. First, these solutions can be represented in non-privileged coordinate systems, and their identification can be a nontrivial task. Second, there is no guarantee that all the corresponding non-equivalent exact solutions of the field equations are known. These remarks are true for both Stackel spaces and Petrov spaces. Therefore, a direct classification based on obtaining all relevant exact solutions is necessary. We emphasize that the main purpose of the classification is not to find new exact solutions of Einstein's equations, but to enumerate all exact solutions possessing a given symmetry.

In papers \cite{33}, \cite{34}, \cite{35}, \cite{36}, \cite{37}, \cite{38} (see also the bibliography presented therein) the classification of the Stackel spaces of vacuum and electrovacuum has been carried out for all types of Stackel spaces. The classification of the Stackel spaces of vacuum and electrovacuum, including the cases with cosmological term, for all types of Stackel spaces, except for the case of Einstein spaces with a non-isotropic complete set consisting of three mutually commuting Killing fields, is carried out. According to the terminology adopted in the theory of complete separation of variables, such spaces are called Stackel spaces of type (3.0) (as already mentioned, these spaces are also Petrov spaces of the first type according to Bianchi's classification). In the present paper, the classification of Einstein-Stackel spaces of this type is carried out. It should be considered as a continuation of the \cite{38}, in which the corresponding classification for the vacuum Einstein-Maxwell equations was obtained.

The problem of classification of spaces with symmetry is divided into three stages.

At the first stage all non-equivalent space-time metrics of these spaces in privileged coordinate systems should be enumerated. For Petrov spaces this is done in \cite{2} (see also \cite{38a}).

At the second stage, in the papers \cite{39}, \cite{40}, \cite{41}, \cite{42}. \cite{43}, \cite{44}, \cite{45}, all non-equivalent potentials of admissible electromagnetic fields, which do not violate the conditions of non-commutative integration, are found and all admissible electromagnetic fields have been classified for Petrov spaces with simply transitive three-parameter groups of motions.

At the third stage the classification of exact solutions of Maxwell, Einstein and Einstein-Maxwell equations should be realized. For the Stackel spaces it is realized in the above papers. For homogeneous Einstein-Petrov spaces with four-parameter group of motions the problem is solved in \cite{46}, \cite{47}, \cite{48}, \cite{49}. In \cite{50}, using numerical methods, solutions of the Einstein-Maxwell equations for Petrov spaces with Bianchi type IX motion group are obtained. For other Petrov spaces, only exact solutions of the source-free Maxwell equations have been classified (see \cite{51}, \cite{52}, \cite{53}, \cite{54}).

Note that all Ricci-flat Stackel spaces of type (3.0), as well as Einstein Stackel spaces of other types, have already been found as special cases of electrovacuum Stackel spaces (see also \cite{56}, \cite{57}, \cite{58}, \cite{59}, \cite{60}, \cite{61}).

\section{ System of Einstein equations}

\noindent
The system of Einstein's vacuum equations with the cosmological term $\Lambda$ has the form:
\begin{equation} \label{1}
R_{ij}=g_{ij} \Lambda,
\end{equation}
where
$$
R_{ir} =\Gamma _{ir,\ell }^{\ell } -\Gamma_{\ell i,k}^{\ell } + \Gamma_{ir}^{\ell} \Gamma _{\ell j}^{j} -\Gamma_{ij}^{\ell } \Gamma_{r\ell }^{j}.
$$
is the Ricci tensor, $g_{ij} $ is the metric tensor of spacetime $V_{4}$.

\noindent
 If the metric of the space $V_{4}$ is invariant with respect to the abelian group of motions $G_3$ given by the triple of Killing vectors:
$$
\xi _{\alpha }^{i} =\delta _{\alpha }^{i},
$$
 then there exists a coordinate system $\left\{u^{i} \right\}$ in which the components of the metric tensor $g_{ij} $ depend only on the variable $u^{0} $. Such a coordinate system is called a privileged coordinate system. In the case where the group $G_3$ acts simply transitively on a non-isotropic hypersurface of the space $V_{4} $, the space is called a non-isotropic Petrov space of the first type according to the Bianchi classification (in the theory of Stackel spaces -- a Stackel space of type (3.0)).
 The metric tensor of a Stackel space of type (3.0) in the canonical privileged coordinate system has the form:
\begin{equation} \label{2}
g_{ij} =-\varepsilon \delta _{i}^{0} \delta _{j}^{0} +\delta _{i}^{\alpha } \delta _{j}^{\beta } g_{\alpha \beta } \left(u^{0} \right)du^{\alpha } du^{\beta } \qquad \qquad \left(\varepsilon =\pm 1\right).
\end{equation}
Here and hereafter, the designations used are:
$$
i,j,k,e=0,1,2,3;\qquad \alpha,\beta,\gamma =1,2,3.
$$
By the repeated upper and lower indices we summarize within the change of indices. The determinant of the matrix $g_{\alpha \beta } $ can be represented as:
\begin{equation} \label{3}
\det | g_{\alpha \beta }| =\varepsilon \ell^{2} \qquad (\ell =\ell (u^{0})).
\end{equation}
When $\varepsilon =1$, the variable $u^{0}$ is temporal. If $\varepsilon =-1,\quad u^{0}$ is a spatial variable.
The variables $\left\{u^{\alpha } \right\}$ of the privileged coordinate system $\left\{u^{i} \right\}$ are defined with exactness to the group of admissible transformations of the form:
\begin{equation} \label{4}
\widetilde{u}^{\alpha } =S_{\beta }^{\alpha } u^{\beta }, \quad S_{\alpha }^{\beta} =const,\quad \det|S_{\alpha }^{\beta}| \ne 0.
\quad \end{equation}
The components of the $R_{j}^{i} $ of the Ricci tensor are given in the book \cite{61}
\begin{equation} \label{5}
R_{0}^{0} =\frac{\varepsilon }{4} \left(2\widetilde{\kappa }_{\alpha ,0}^{\alpha } +\widetilde{\kappa }_{\beta }^{\alpha } \widetilde{\kappa }_{\alpha }^{\beta } \right),\qquad R_{\beta }^{\alpha } =\frac{\varepsilon }{2\ell } \left(\ell \widetilde{\kappa }_{\alpha }^{\beta } \right)_{,0} ,
\end{equation}
where \quad $\widetilde{\kappa }_{\beta }^{\alpha } =g^{\alpha \gamma } g_{\gamma \beta ,0}. $

\noindent
Let us perform the following substitution of the variable $u^{0}$:
\begin{equation} \label{6}
d u^{0} (\tau )/d\tau=\ell.
\end{equation}
We will denote the derivatives on the variable $\tau $ by points. Let's denote:
\begin{equation} \label{7}
\kappa _{\beta }^{\alpha } =g^{\alpha \gamma } \dot{g}_{\gamma \beta }
\end{equation} and introduce a $\gamma$ function that satisfies the condition:
\begin{equation} \label{8}
\ddot{\gamma }=\xi \ell^2\exp{\lambda} \quad (2\Lambda =\varepsilon \xi \exp{\lambda}).
\end{equation}
Then Einstein's \eqref{1} equations can be represented as:
\begin{equation} \label{9}
 \dot{\kappa}_{\beta}^{\alpha } =\delta_{\beta}^{\alpha}\ddot{\gamma},
\end{equation}
\begin{equation} \label{10}
 2\dot{\kappa }_{\alpha }^{\alpha } +\kappa _{\beta }^{\alpha } \kappa _{\alpha }^{\beta } -(\kappa _{\alpha }^{\alpha })^{2}
=2\ddot{\gamma}. \end{equation}
The system \eqref{9} can be integrated. After the first integration we obtain:
\begin{equation} \label{11}
\kappa _{\beta }^{\alpha } =\delta _{\beta}^{\alpha }\dot{\gamma }+C_{\beta }^{\alpha } \qquad \left(C_{\beta }^{\alpha } =const\right)
\end{equation}
Using the \eqref{7} relation, we obtain a system of equations:
\begin{equation} \label{12}
\dot{g}_{\alpha \beta } =g_{\beta}^{\alpha }\dot{\gamma }+C_{\beta }^{\alpha } g_{\alpha \gamma }
\end{equation}
Let's denote:
$$
g_{\alpha \beta } =\eta_{\alpha \beta } \exp{\gamma}
$$
Then the system \eqref{11} can be represented as:
\begin{equation} \label{13}
\dot{\eta}_{\alpha \beta } =C_{\beta }^{\alpha } \eta_{\alpha \gamma }.
\end{equation}
Since the $\gamma$ function is not included in the \eqref{13} equations, the system can be called the first autonomous subsystem.

\quad

\qquad

\section{Solution of the first autonomous subsystem}

The matrix $C^\alpha_\beta$ can be simplified using the following transformations of the coordinate system $\{\hat{u}\},$ acting on the hypersurface of transitivity $V_3\in V_4$:
\begin{equation}\label{14}
\tilde{u}^\alpha = S^\alpha_\beta u^\beta, \quad S^\alpha_\beta =const, \quad \det{\hat{S}} \ne 0.
\end{equation}
In the new coordinate system, the matrix $\tilde{C^\alpha_\beta}$ will have the form:
\begin{equation}\label{15}
\tilde{C^\alpha_\beta} = S^\alpha_\gamma C^\gamma_\rho (\hat{S}^{-1})^\rho_\gamma.
\end{equation}
In the general case, the matrix $C^\alpha_\beta$ is not equivalent to a diagonal matrix with respect of the transformations \eqref{15}. However, it can be shown that there exist coordinate transformations of the form \eqref{14} which allow to zeroize the non-diagonal elements of the matrix $C^\alpha_\beta$ corresponding to any given index $\alpha$, e.g. $\alpha =1$:
\begin{equation}\label{16}
C^1_p = C^q_1 = 0 \quad (p,q =2,3).
\end{equation}
Let us assume that such a transformation is realized and the matrix $C^\alpha_\beta$ has the form:
\begin{equation}\label{17}
C^\alpha_\beta = C^1_1 \delta^\alpha_1\delta^1_\beta + C_p^q \delta^p_\beta \delta_q^\alpha.
\end{equation}
We substitute the matrix \eqref{17} into the system of equations \eqref{13}. As result these equations take the form:
\begin{equation}\label{18}
\dot{\eta}_{11} =C_1^{1}\eta_{11}, \quad
\dot{\eta}_{pq} =C_p^{p'}\eta_{p'q}.
\end{equation}
By virtue of the symmetry of the matrix $\eta_{pq}$, from the second equation of the system \eqref{18} it follows:
\begin{equation}\label{19}
(C_2^{2} - C_3^{3})\eta_{23} + C_3^{2}\eta_{33} - C_2^{3}\eta_{22}=0.
\end{equation}
So if
\begin{equation}\label{20}
(C_2^{2} - C_3^{3})^2 + (C_3^{2})^2 + (C_2^{3})^2 \ne 0,
\end{equation}
the functions $\eta_{pq}$ are linearly dependent with constant coefficients. It is not difficult to show that from the equation\eqref{18} their linear dependence follows also when the condition \eqref{20} is satisfied. Therefore, in order to realize a complete classification of solutions to equation\eqref{19}, we need to consider two variants:

\noindent
- the linear function space $\eta_{pq}$ is one-dimensional,

\noindent
- the linear space of functions $\eta_{pq}$ is two-dimensional.

Considering these variants, we obtain three non-equivalent solutions of the equation \eqref{19}:

\quad

\noindent 1.
\begin{equation} \label{21}
\eta_{pq} =\delta_{pq}\eta_{pp}, \quad C_3^{2}= C_2^{3} =0.
\end{equation}

\quad

\noindent 2.
\begin{equation} \label{22}
\eta_{33} = -\eta_{22}, \quad C_2^{2} = C_3^{3}, \quad C_3^{2}=-C_2^{3}.
\end{equation}

\quad

\noindent 3.
\begin{equation} \label{23}
\eta_{33} = 0, \quad C_2^{2} = C_3^{3}, \quad C_3^{2}=0.
\end{equation}

\quad

Let us use the obtained classification to integrate the system of equations\eqref{18}.

\quad

{\noindent}
{\bf A.}\quad Variant \eqref{21}. We substitute \eqref{21} into the system \eqref{18}. This results in the following solution:
\begin{equation} \label{24}
\eta_{11} =\varepsilon_1\exp \left(C_{1}^{1} \tau \right),\qquad \eta _{22} =\varepsilon_2\exp (C_{2}^{2}\tau),\qquad \eta _{33} =\varepsilon_3\exp (C_{3}^{3}\tau).
\end{equation}
The matrix $C^\alpha_\beta$ has the form:
$$
C^{\alpha}_{\beta} =\left(\begin{matrix}& {c_1} & {0} & {0} \\ &{0} & {c_2} & {0} \\ &{0} & {0} & {c_3} \end{matrix}\right) =const.
$$

{\noindent}
{\bf B.}\quad Variant \eqref{22}. Let us substitute \eqref{21} into the system \eqref{18}. This results in the following system:
\begin{equation}\label{25}
\dot{\eta}_{22}=c\eta_{22}+b\eta_{23}, \quad \dot{\eta}_{22}=-b\eta_{22}+c\eta_{23}, \quad \eta_{33}=-\eta_{22},
\end{equation}
where \quad $c=C^2_2, \quad b=C^3_2$.
The solution can be represented as:
\begin{equation}\label{26}
\eta_{11}=\exp(\sigma+c), \quad \eta_{22}=\exp\sigma\sin\omega, \quad \eta_{23}=\exp\sigma\cos\omega,
\end{equation}
here $\sigma = \tau\cos{a}, \quad \omega = \tau\sin{a}$;
\begin{equation} \label{27}
C^{\alpha}_{\beta} =\left(\begin{matrix}& {c+\cos{a}} & {0} & {0} \\ &{0} & {\cos{a}} & {-\sin{a}} \\ &{0} & {\sin{a}} & {\cos{a}} \end{matrix}\right).
\end{equation}

\quad

{\noindent}
{\bf C.}\quad Variant \eqref{23}.

The system of equations \eqref{18} takes the form:
\[\left\{\begin{array}{c} {\dot{\eta}_{22} =a \eta_{22} + b\eta_{23}, \qquad \qquad \qquad \qquad \qquad } \\\\{\dot{\eta}_{23} =a^{2} \eta_{23} \qquad \qquad \qquad \qquad \qquad \qquad \qquad } \end{array}\right. \]
The result of the integration is:

\noindent $\eta_{23} =p\exp \left(a^{2} \tau \right),\qquad \eta_{22} =q\exp \left(a^{2} \tau \right)+\tau a^{3} \eta_{23} $, which is equivalent to the solution:
\begin{equation} \label{28}
\eta_{23} =\exp \left(\tau \right),\qquad \eta_{22} =\tau \exp \tau \qquad \eta_{33}=0.
\end{equation}

The matrix $C^\alpha_\beta$ has the form:
$$
C^{\alpha}_{\beta} =\left(\begin{matrix}& {c} & {0} & {0} \\ &{0} & {1} & {0} \\ &{0} & {1} & {1} \end{matrix}\right).
$$
The matrix $C^\alpha_\beta$ has the form:
$$
C^{\alpha}_{\beta} =\left(\begin{matrix}& {c} & {0} & {0} \\ &{0} & {1} & {0} \\ &{0} & {1} & {1} \end{matrix}\right).
$$

\section{Solution of the second autonomous subsystem}

 Let us consider the remaining equations from the system of Einstein equations \eqref{1}. From equation \eqref{11}:
\begin{equation}\label{29}
\ell^2=\exp(3\gamma + \tau C^\alpha_\alpha +k) \quad (k=const).
\end{equation}
Now the equation \eqref{8} can be represented as:
\begin{equation}\label{30}
\ddot{\gamma} = \xi\exp(3\gamma + \tau C^\alpha_\alpha + \lambda +k)
\end{equation}
Substituting \eqref{11} into equation \eqref{12}, we get:
\begin{equation} \label{31}
4\ddot{\gamma} -6{\dot{\gamma}}^2 -4\dot{\gamma}C^\alpha_\alpha +C^\alpha_\beta C^\beta_\alpha - (C^\alpha_\alpha )^2 =0.
\end{equation}
Let's denote:
\begin{equation}\label{33}
\varphi =3\gamma + \tau C^\alpha_\alpha + \lambda +k, \quad (C^\alpha_\alpha )^2-3C^\alpha_\beta C^\beta_\alpha =2\epsilon p^2, \quad \epsilon =0, \pm1, \quad p\ne0.
\end{equation}
This results in two equations on the same function $\varphi$, which include the components of the constant matrix $C_\beta^\alpha$:
\begin{equation}\label{32}
\left\{\begin{array}{c} {4\ddot{\varphi} -2{\dot{\varphi}}^2 +3C^\alpha_\beta C^\beta_\alpha - (C^\alpha_\alpha )^2 =0, \qquad \qquad \qquad \qquad } \\ {\ddot{\varphi}=3\xi\exp{\varphi}. \qquad \qquad \qquad \qquad \qquad \qquad \qquad } \end{array}\right. \end{equation}
From the first equation of the \eqref{33} system:
\begin{equation}\label{34}
\int\frac{d(\dot{\varphi})}{(\dot{\varphi}^2 + \epsilon p^2)}=\frac{\tau}{2}.
\end{equation}
Depending on the value of the parameter $\epsilon $, the solutions of the equations \eqref{32} are of the form:

\quad

\noindent 1)\quad $\epsilon=1$
\begin{equation}\label{35}
\exp{\varphi}=\frac{p^2}{6\cos^2(\frac{\tau p}{2})}, \quad \xi=1;
\end{equation}
$$
g_{\alpha\beta}=\eta_{\alpha\beta}\left( \frac{p^2}{6\cos^2(\frac{\tau p}{2})}\right)^{\frac{1}{3}}\exp{-\left( \frac{\lambda +c\tau+3\sigma}{3}\right)}.
$$
\noindent 2)\quad $\epsilon=-1$.
\begin{equation}\label{36}
\exp{\varphi}=\frac{p^2}{6\cosh^2(\frac{\tau p}{2})}, \quad \xi=-1;
\end{equation}
$$
g_{\alpha\beta}=\eta_{\alpha\beta}\left( \frac{p^2}{6\cosh^2(\frac{\tau p}{2})}\right)^{\frac{1}{3}}\exp{-\left( \frac{\lambda +c\tau+3\sigma}{3}\right)}
$$

\noindent 3)\quad $\epsilon=-1$.
\begin{equation}\label{37}
\exp{\varphi}=\frac{p^2}{6\sinh^2(\frac{\tau p}{2})}, \quad \xi=1;
\end{equation}
$$
g_{\alpha\beta}=\eta_{\alpha\beta}\left( \frac{p^2}{6\sinh^2(\frac{\tau p}{2})}\right)^{\frac{1}{3}}\exp{-\left( \frac{\lambda +c\tau+3\sigma}{3}\right)}.
$$

\noindent 4)\quad $\epsilon=0$.
\begin{equation}\label{38}
\exp\varphi=\frac{2}{3\tau^2}, \quad \xi=1;
\end{equation}
$$
g_{\alpha\beta}=\eta_{\alpha\beta}\left(\frac{2}{3{\tau}^2}\right)^{\frac{1}{3}}\exp{-\left(\frac{\lambda +c\tau+3\sigma}{3}\right)}.
$$
The components of the matrices $\eta_{\alpha\beta}$ and $C_{\alpha}^{\beta}$ are represented by the expressions \eqref{21}-\eqref{28}. Using these expressions, we find all metrics of the sought Einstein spaces.

\section{Metrics of Einstein spaces}

\noindent
{\bf Variant A.}\quad The matrices $\eta_{\alpha\beta}$ and $C_{\alpha}^{\beta}$ are represented by the expressions \eqref{24}. Let's represent the expression \eqref{33} in the form:
\begin{equation}\label{41}
2((a_1 + a_2 + a_3)^2 - 3(a_1^2 + a_2^2 + a_3^2)) = (2a_1 - a_2 - a_3)^2 +3((a_2 - a_3)^2)) {\Rightarrow}
\end{equation}
$$(2a_1 - a_2 - a_3)^2 +3(a_2 - a_3)^2=- 4\epsilon p^2.$$
With this in mind, we refine the expressions \eqref{35} - \eqref{28}. As a result, we obtain the following metrics $$\left( c= (a_1 + a_2 + a_3), \quad p = \sqrt{{a_1}^2 + {a_1}^2 + {a_1}^2 +a_1(a_2 + a_3)+a_2a_3} \quad \right):$$

\quad

\noindent

\noindent $\bf{1.}\quad\bf{\epsilon=\xi=-1.} $

\begin{equation}\label{40}
ds^2=-\frac{{e_1}{e_2}{e_3} p^2}{6\exp\lambda\cosh^2(\frac{\tau p}{2})}  d\tau^2  +\left(\frac{p^2}{6\cosh^2(\frac{\tau p}{2})\exp(\lambda +c\tau)}\right)^{\frac{1}{3}}(e_\alpha\exp{\tau a_\alpha}du^\alpha du^\alpha),
\end{equation}

\quad

\noindent $\bf{2.}\quad\bf{\epsilon=-1,\quad \xi=1.} $

\begin{equation}\label{40a}
ds^2=-\frac{{e_1}{e_2}{e_3} p^2}{6\exp\lambda\sinh^2(\frac{\tau p}{2})}  d\tau^2  +\left(\frac{p^2}{6\sinh^2(\frac{\tau p}{2})\exp(\lambda +c\tau)}\right)^{\frac{1}{3}}(e_\alpha\exp{\tau a_\alpha}du^\alpha du^\alpha),
\end{equation}

\quad

\noindent $\bf{3.}\quad\bf{\epsilon=-1,\quad \xi=1.} $
\begin{equation}\label{40}
ds^2=-\frac{2{e_1}{e_2}{e_3} d\tau^2}{3{\tau}^{2}\exp\lambda} +\left(\frac{ 2{e_1}{e_2}{e_3} }{3{\tau}^{2}\exp\lambda}\right)^{\frac{1}{3}}(e_\alpha {du^\alpha}^2).
\end{equation}

\quad

\noindent
{\bf Variant B.}\quad The matrices $\eta_{\alpha\beta}$ and $C_{\alpha}^{\beta}$ are represented by the expressions \eqref{26} - \eqref{27}.

\noindent In this case, the expression \eqref{33} takes the form:
\begin{equation}\label{39}
3\cos{\sigma}-c^2 = \epsilon p^2
\end{equation}
Let's refine the expressions \eqref{35} - \eqref{28}. As a result, we obtain the following metrics.

\quad

\noindent $\bf{4.}\quad\bf{\epsilon=\xi=1.} $

\begin{equation}\label{40}
ds^2=\frac{p^2}{6\exp\lambda\cos^2(\frac{\tau p}{2})}  d\tau^2  +\left(\frac{p^2}{6\cos^2(\frac{\tau p}{2})\exp(\lambda +c\tau)}\right)^{\frac{1}{3}}(\exp(c\tau){du^1}^2
\end{equation}
$$
+ \sin\omega({du^2}^2-{du^3}^2) + 2\cos\omega du^2du^3) \quad (p = \sqrt{(3\cos{\sigma}-c^2)}).
$$

\noindent $\bf{5.}\quad\bf{\epsilon=\xi=-1.} $
\begin{equation}\label{40}
ds^2=\frac{p^2}{6\exp\lambda\cosh^2(\frac{\tau p}{2})}  d\tau^2  +\left(\frac{p^2}{6\cosh^2(\frac{\tau p}{2})\exp(\lambda +c\tau)}\right)^{\frac{1}{3}}(\exp(c\tau){du^1}^2
\end{equation}

\noindent $\bf{6.}\quad\bf{\epsilon=-\xi=-1.} $
\begin{equation}\label{40}
ds^2=\frac{p^2}{6\exp\lambda\sinh^2(\frac{\tau p}{2})}  d\tau^2  +\left(\frac{p^2}{6\sinh^2(\frac{\tau p}{2})\exp{(\lambda +c\tau)}}\right)^{\frac{1}{3}}(\exp(c\tau){du^1}^2
\end{equation}
$$
+ \sin\omega({du^2}^2-{du^3}^2) + 2\cos\omega du^2du^3) \quad (p = \sqrt{(c^2 - 3\cos{\sigma})}).
$$

\noindent $\bf{7.}\quad\bf{\epsilon=0,\quad \xi=1.} $
\begin{equation}\label{41}
ds^2=\frac{2}{3\tau^2\exp\lambda}d\tau^2 +\left(\frac{2}{3\tau^2\exp(\lambda+c\tau)}\right)^{\frac{1}{3}}(\exp c\tau{du^1}^2
\end{equation}
$$
+ \sin\omega({du^2}^2-{du^3}^2) + 2\cos\omega du^2du^3) \quad (c = \sqrt{3}\sin{a}).
$$

\noindent
{\bf Variant C.}\quad The matrices $\eta_{\alpha\beta}$ and $C_{\alpha}^{\beta}$ are represented by the expressions \eqref{28}. Let's represent the \eqref{33} expression as:
$$
b^2=-\epsilon p^2
$$

and refine the expressions \eqref{35} - \eqref{28}. As a result, we obtain the following metrics.

\quad\noindent ${\bf{8.}}\quad\epsilon=-\xi=1,\quad p=b. $

\begin{equation}\label{40}
ds^2=-\frac{{e_1}{e_2}{e_3} p^2 d\tau^2}{6\exp(\lambda)\cosh^2(\frac{\tau p}{2})} +\left(\frac{p^2}{6\cosh^2(\frac{\tau p}{2})\exp(\lambda -\tau)}\right)^{\frac{1}{3}}({du^1}^2\exp{-p\tau}
\end{equation}
$$+c\tau {du^2}^2 + 2du^2du^3)$$,
\quad

\noindent ${\bf9.}\quad \epsilon=-1,\quad \xi=1,\quad p=b. $.
\begin{equation}\label{40}
ds^2=-\frac{{e_1}{e_2}{e_3} p^2 d\tau^2}{6\exp(\lambda)\sinh^2(\frac{\tau p}{2})} +\left(\frac{p^2}{6\sinh^2(\frac{\tau p}{2})\exp(\lambda -\tau)}\right)^{\frac{1}{3}}({du^1}^2\exp{-p\tau}
\end{equation}
$$+c\tau {du^2}^2 + 2du^2du^3)$$,

\quad

\noindent ${\bf10.}\quad \epsilon=-1,\quad \xi=1,\quad b=0. $.
\begin{equation}\label{40}
ds^2=\frac{2 d\tau^2}{3{\tau}^{2}\exp\lambda} +\left(\frac{2}{3{\tau}^{2}\exp\lambda}\right)^{\frac{1}{3}}({du^1}^2+c\tau {du^2}^2 + 2du^2du^3).
\end{equation}

\quad

\section{Conclusion}
 In this paper we consider the problem of classification of Einstein spaces for the case when the three-parameter abelian group of motions acts on non-null hypersurfaces. All non-equivalent metrics with respect to the group of admissible coordinate transformations are enumerated. Thus, the classification of Einstein spaces in which the geodesic equations can be integrated by the method of complete separation of variables in the Hamilton-Jacobi equation is completed. Since flat Ricci spaces of this type were classified already in the middle of the 20th century (see, e.g., \cite{61}), the problem of classifying exact solutions of the vacuum Einstein equations for all Stackel spaces is completely solved.

\quad

\end{document}